\documentclass[11pt]{article}

\usepackage[numbers,sort&compress]{natbib}
\usepackage{color}
\usepackage{amsfonts}
\usepackage{mathrsfs}
\usepackage{bm}
\usepackage{amssymb,amsmath} 
\usepackage{amsthm}

\allowdisplaybreaks

\catcode`!=11
\let\!int\int \def\int{\displaystyle\!int}
\let\!lim\lim \def\lim{\displaystyle\!lim}
\let\!sum\sum \def\sum{\displaystyle\!sum}
\let\!sup\sup \def\sup{\displaystyle\!sup}
\let\!inf\inf \def\inf{\displaystyle\!inf}
\let\!cap\cap \def\cap{\displaystyle\!cap}
\let\!max\max \def\max{\displaystyle\!max}
\let\!min\min \def\min{\displaystyle\!min}
\let\!frac\frac \def\frac{\displaystyle\!frac}
\catcode`!=12

\let\oldsection\section

\renewcommand\section{\setcounter{equation}{0}\oldsection}

\allowdisplaybreaks

\def\Li{\operatorname{Li}}

\def\pf{\noindent{\bf Proof.}\rm\quad}

\def\N{\mathbb{N}}

\newtheorem{thm}{Theorem}[section]
\newtheorem{lem}[thm]{Lemma}
\newtheorem{cor}[thm]{Corollary}
\newtheorem{re}{Remark}[section]

\begin{document}


\title {\bf On $q$-analogues of quadratic Euler sums}

\author{
{Zhonghua Li\thanks{zhonghua\_li@tongji.edu.cn (Z. Li) }\quad Ce Xu\thanks{Corresponding author. Email: xuce1242063253@163.com (C. XU)}}\\[1mm]
\small $*$ School of Mathematical Sciences, Tongji University\\
\small Shanghai 200092, P.R. China\\
\small$\dagger$ School of Mathematical Sciences, Xiamen University\\
\small Xiamen
361005, P.R. China}

\date{}

\maketitle \noindent{\bf Abstract } In this paper we study $q$-analogues of Euler sums and present a new family of identities by using the method of Jackson $q$-integral representations of series. We then apply it to obtain a family of identities relating quadratic Euler sums to linear sums and $q$-polylogarithms. Furthermore, we also use certain stuffle products to evaluate several $q$-series with $q$-harmonic numbers. Some interesting new results and illustrative examples are considered. Finally, let $q$ tend to $1$, we obtain some explicit relations for the classical Euler sums.
\\[2mm]

\\[2mm]

\tableofcontents
\section{Introduction}

For positive integers $m$ and $k$, let $H^{(k)}_m$ and ${\overline{H}}^{(k)}_m$ stand for the $m$-th generalized harmonic number and the $m$-th generalized alternating harmonic number defined by \cite{A1985,FS1998,X2016}
\begin{align*}
H_m^{(k)}: = \sum\limits_{j=1}^m {\frac{1}{{{j^k}}}} ,\qquad\qquad \overline{H}_m^{(k)}: = \sum\limits_{j=1}^m {\frac{(-1)^{j-1}}{{{j^k}}}},
\end{align*}
respectively. If $k>1$, the generalized harmonic number $H^{(k)}_m$ converges to the (Riemann) zeta value $\zeta(k)$:
$$\lim\limits_{m\rightarrow \infty}H^{(k)}_m=\zeta(k).$$
When $k=1$, $H_m^{(1)}\equiv H_m$ (resp. $\overline{H}_m^{(1)}\equiv \overline{H}_m$) is the classical harmonic number (resp. the classical alternating harmonic number).

Let $n$ be a positive integer. Let $k_1,\ldots,k_n$ be nonzero integers and let $k$ be a positive integer with $k\geq 2$. The classical Euler sums are defined by the convergent series
\begin{align*}
&S(k_1,\ldots,k_n;k): = \sum\limits_{m=1}^\infty\frac{X_m(k_1)\cdots X_m(k_n)}{m^k},\\
&\overline{S}(k_1,\ldots,k_n;k):=\sum\limits_{m= 1}^\infty\frac{X_m(k_1)\cdots X_m(k_n)}{m^k}(-1)^{m-1},
\end{align*}
where
$$X_m(k): = \begin{cases}
   H_m^{(k)} & \text{if\;} k \geq 1,  \\
   \overline{H}_m^{(-k)} & \text{if\;}k \leq -1.
\end{cases}$$
Here we call $|k_1|+\cdots+|k_n|+k$ the weight, and $n$ the depth. Throughout the paper,  for a positive integer $k$, we use ${\bar k}$ to denote the negative entry $-k$. For example, we have
\begin{align*}
&S(1,\bar{2},3;4) = S(1,-2,3;4),\qquad \qquad\overline{S}(\bar{1},\bar{2},3;4)=\overline{S}(-1,-2,3;4).
\end{align*}

It is clear that every Euler sum of weight $w$ and depth $n$ is clearly a $\mathbb{Q}$-linear combination of multiple zeta values or multiple zeta star values (that is, values of multiple zeta functions or multiple zeta star functions at integer arguments) of weight $w$ and depth less than or equal to $n+1$. In other words, multiple zeta (star) values are ``atomic" quantities into which Euler sums decompose. The multiple zeta and zeta star values are defined by (\cite{H1997,KP2013,Z2014})
\begin{align*}
&\zeta( \mathbf{k})\equiv\zeta(k_1, \ldots, k_n):=\sum\limits_{m_1>\cdots>m_n\geq 1}\prod\limits_{j=1}^n m_j^{-|k_j|}{\rm sgn}(k_j)^{m_j},\\
&\zeta^\star(\mathbf{k})\equiv\zeta^\star(k_1, \ldots, k_n):=\sum\limits_{m_1\ge\cdots  \ge m_n \ge 1} \prod\limits_{j=1}^n m_j^{-|k_j|}{\rm sgn}(k_j)^{m_j},
\end{align*}
where for convergence $|k_1|+\cdots+|k_j|> j$ for $j= 1, 2, \ldots, n$, and
\[{\rm sgn}(k_j):=\begin{cases}
   1  & \text{\;if\;} k_j>0,  \\
   -1, & \text{\;if\;} k_j<0.
\end{cases}\]
Here, we call $l(\mathbf{k}):=n$ and $|\mathbf{k}|:=\sum\limits_{j=1}^n|k_j|$ the depth and the weight of multiple zeta values, respectively.

Euler sums and multiple zeta values have connections with many branches of mathematics; see especially Zagier \cite{DZ1994}. The evaluation of Euler sums also has been useful in various areas of theoretical physics, including in support of Feynman diagram calculations and in resolving open questions on Feynman diagram contributions and relations among special
functions \cite{C2005,C2008}, including the dilogarithm, Clausen function, and generalized hypergeometric function. An array of harmonic number sums and multiple zeta values is required in calculations
of high energy physics. These quantities appear for instance in developing the scattering theory of massless quantum electrodynamics \cite{Bl1999}.
Broadhurst (see Borwein and Girgensohn \cite{BG1996}) encountered them in relation with Feynman diagrams and associated knots in perturbative quantum field theory.

A good deal of work on Euler sums has been focused on the problem of determining when complicated sums can be expressed in terms of simpler sums. Thus, researchers are interested in determining which sums can be expressed in terms of other sums of lesser depth.
The origin of the study of Euler sums went back to the correspondence of Euler with Goldbach in 1742-1743 (see \cite{H2007}) and Euler's paper \cite{E1775} that appeared in 1776. Euler studied linear (or double) Euler sums and established some important formulas for them. For example, he proved that (see \cite{BBG1994,FS1998})
$$S(1;k)=\frac{1}{2}\left\{ (k + 2)\zeta(k + 1) - \sum\limits_{i = 1}^{k - 2}\zeta(k - i)\zeta(i + 1)\right\}.$$
Moreover, Euler proved that the linear sums $S(l;k)$ ($l\geq1, k\geq2$) are reducible to zeta values whenever $k+l$ is less that $7$ or when $k+l$ is odd and less than $13$. Furthermore, he conjectured that the linear sums $S(l;k)$ would be reducible to zeta values whenever $k+l$ is odd, and even proposed the general formula. In \cite{BBG1995}, D. Borwein, J. M. Borwein and R. Girgensohn proved the conjecture,
and in \cite{BBG1994}, D. H. Bailey, J. M. Borwein and R. Girgensohn demonstrated that it is \lq\lq very likely" that the linear sums $S(l;k)$ with $k+l> 7$ and $k+l$ even, are not reducible. After that many different methods, including partial fraction expansions, Eulerian Beta integrals, summation formulas for generalized hypergeometric functions and contour integrals, have been used to evaluate these sums (see \cite{BBG1994,BBG1995,FS1998}). For example, P. Flajolet and B. Salvy informed us about some ongoing work of theirs (\cite{FS1998}) to evaluate Euler sums in an entirely different way, namely using contour integration and the residue theorem. In this way they manage to prove, for example, that the sums $S(1,1,1;k)$ with $k=2,3,4,6$ can be evaluated in terms of zeta values. There are also a lot of recent contributions on nonlinear Euler sums (depth $\geq 2$), see \cite{X2017,X2016}. For example, in \cite{X2017}, we proved that all Euler sums of the form $S(k_1,k_2;k)$ with weight $4,5,6,7,9$ are expressible polynomially
in terms of zeta values. For weight $8$, all such sums are the sum of a polynomial in zeta values and a rational multiple of $S(2;6)$. And all weight $10$ quadratic sums $S(1,l;k)$ are reducible to $S(2;6)$ and $S(2;8)$.

So far, surprising little work has been done on $q$-analogues of Euler sums and multiple zeta values. Actually, there are many possible ways to $q$-extend the Euler sums and multiple zeta values. Here we recall one $q$-analogue. Let $q$ be a fixed real number with $0<q<1$. Let $n$ be a positive integer. For a sequence ${\mathbf k}=(k_1,...,k_n)$ of positive integers, a sequence $\mathbf{x}=(x_1,\ldots,x_n)$ of variables with $-1\leq x_i\leq1$, a positive integer $k$ and a variable $x$ with $-1<x<1$,  we set
\begin{align}
S\left[ {\left. {\begin{array}{*{20}{c}}
   {\bf{k}}  \\
   {\bf{x}}  \\
\end{array}} \right|\begin{array}{*{20}{c}}
   k  \\
   x  \\
\end{array}} \right] \equiv S\left[ {\left. {\begin{array}{*{20}{c}}
   {{k_1}, \ldots ,{k_n}}  \\
   {{x_1}, \ldots ,{x_n}}  \\
\end{array}} \right|\begin{array}{*{20}{c}}
   k  \\
   x  \\
\end{array}} \right] := \sum\limits_{m = 1}^\infty  {\frac{{{\zeta _m}[{k_1},{x_1}] \cdots {\zeta _m}[{k_n},{x_n}]}}{{{{[m]}^k}}}} {x^m},
\label{Eq:q-EulerSum}
\end{align}
where $[m]$ denotes the $q$-analogue of a nonnegative integer, defined by
$$[m]\equiv[m]_q:=\frac{1-q^m}{1-q},$$
 and ${\zeta _m}[k,x]$ is the partial sum of the $q$-polylogarithm function ${\rm Li}_k[x]$, defined as
$$\zeta_m[k,x]:=\sum\limits_{j=1}^m\frac{x^j}{[j]^k}.$$
Here the $q$-polylogarithm function ${\rm Li}_k[x]$ is defined by
$${\rm Li}_k[x]:=\sum\limits_{m=1}^{\infty} \frac{x^m}{[m]^k},\quad (-1<x<1).$$
Note that
$$\ln[1-x]:=-\Li_1[x]$$
is the $q$-analogues of natural logarithm function. If $n=0$ in \eqref{Eq:q-EulerSum}, we set
$$S\left[\left.\begin{matrix}
\varnothing\\
\varnothing
\end{matrix}\right|\begin{matrix}
k\\
x
\end{matrix}\right]:={\rm Li}_k[x].$$
When taking the limit $q\rightarrow 1$ and $x\rightarrow 1$ with $x_j=1$ in \eqref{Eq:q-EulerSum} we get
$$\mathop {\lim }\limits_{q \to 1} S\left[ {\left. {\begin{array}{*{20}{c}}
   {{k_1}, \ldots ,{k_n}}  \\
   {1, \ldots ,1}  \\
\end{array}} \right|\begin{array}{*{20}{c}}
   k  \\
   1  \\
\end{array}} \right] = S\left( {{k_1}, \ldots ,{k_n};k} \right).$$
For a real number $a\neq -1,-2,\ldots$, we set
$$H_k[x,a]: = \sum\limits_{m= 1}^\infty\frac{x^{m + a}}{[m + a]^k},\quad k\in \N,x\in(-1,1),$$
where for a general real number $b$,
\[\left[ b \right] \equiv {\left[ b \right]_q}: = \frac{{1 - {q^b}}}{{1 - q}}.\]

There are fewer results for sums of the type \eqref{Eq:q-EulerSum}. Some related results for $q$-Euler type sums may be seen in the works of \cite{DPP2012,KHP2015,L2015,Sa2014,M2014,W2014,Xu2016,Z2014} and references therein.
The second author jointly with M. Zhang and W. Zhu \cite{Xu2016} proved that for positive integer $k\geq 2$, the q-linear sum
\[S\left[ {\left. {\begin{array}{*{20}{c}}
   1  \\
   1  \\
\end{array}} \right|\begin{array}{*{20}{c}}
   k  \\
   q  \\
\end{array}} \right]\]
can be expressed as a rational linear combination of products of $q$-polylogarithms,
the quadratic sum
\[S\left[ {\left. {\begin{array}{*{20}{c}}
   {1,1}  \\
   {1,1}  \\
\end{array}} \right|\begin{array}{*{20}{c}}
   k  \\
   q  \\
\end{array}} \right]\]
and the cubic combination sum
\[S\left[ {\left. {\begin{array}{*{20}{c}}
   {1,1,1}  \\
   {1,1,1}  \\
\end{array}} \right|\begin{array}{*{20}{c}}
   k  \\
   q  \\
\end{array}} \right] - 3S\left[ {\left. {\begin{array}{*{20}{c}}
   {1,2}  \\
   {1,1}  \\
\end{array}} \right|\begin{array}{*{20}{c}}
   k  \\
   q  \\
\end{array}} \right]\]
are reducible to linear $q$-sums and to polynomials in $q$-polylogarithms. Some simple examples are
\begin{align*}
&S\left[ {\left. {\begin{array}{*{20}{c}}
   1  \\
   1  \\
\end{array}} \right|\begin{array}{*{20}{c}}
   2  \\
   q  \\
\end{array}} \right] = {\rm{L}}{{\rm{i}}_3}\left[ q \right] + {\rm{L}}{{\rm{i}}_3}\left[ {{q^2}} \right],\\
&S\left[ {\left. {\begin{array}{*{20}{c}}
   1  \\
   1  \\
\end{array}} \right|\begin{array}{*{20}{c}}
   3  \\
   q  \\
\end{array}} \right] = \frac{3}{2}{\rm{L}}{{\rm{i}}_4}\left[ {{q^2}} \right] + {\rm{L}}{{\rm{i}}_4}\left[ q \right] - \frac{1}{2}{\rm{Li}}_2^2\left[ q \right],\\
&S\left[ {\left. {\begin{array}{*{20}{c}}
   {1,1}  \\
   {1,1}  \\
\end{array}} \right|\begin{array}{*{20}{c}}
   2  \\
   q  \\
\end{array}} \right] = \frac{7}{2}{\rm{L}}{{\rm{i}}_4}\left[ {{q^2}} \right] + 2{\rm{L}}{{\rm{i}}_4}\left[ q \right] - \frac{1}{2}{\rm{Li}}_2^2\left[ q \right] - \left( {1 - q} \right)\left( {{\rm{L}}{{\rm{i}}_3}\left[ {{q^2}} \right] + {\rm{L}}{{\rm{i}}_3}\left[ q \right]} \right).
\end{align*}

Similarly, the multiple zeta (star) values also exist many $q$-extension. For example, one definition of $q$-multiple zeta (star) values is
\begin{align*}
&{\zeta}\left[ k_1, k_2, \ldots, k_n \right]: = \sum\limits_{m_1>m_2>\cdots>m_n\geqslant 1} {\frac{{{q^{{m_1} + {m_2} +  \cdots  + {m_n}}}}}
{{\left[ {{m_1}} \right]_q^{{k_1}}\left[ {{m_2}} \right]_q^{{k_2}} \cdots \left[ {{m_n}} \right]_q^{{k_n}}}}},\quad  \left(n, k_i\in \N \right),\\
&\zeta^ \star \left[ k_1, k_2, \ldots, k_n \right]: = \sum\limits_{{m_1} \geqslant {m_2} >  \cdots  \geqslant {m_n} \geqslant 1} {\frac{{{q^{{m_1} + {m_2} +  \cdots  + {m_n}}}}}
{{\left[ {{m_1}} \right]_q^{{k_1}}\left[ {{m_2}} \right]_q^{{k_2}} \cdots \left[ {{m_n}} \right]_q^{{k_n}}}}}, \quad \left( n,k_i\in \N \right).
\end{align*}

The other one particularly well-behaved $q$-analog of the multiple zeta functions is defined in \cite{Z2014} by Zhao, generalizing the Riemann $q$-zeta function studied by Kaneko et al. \cite{K2003}. It is very important to understand the relations between their special values, see \cite{Br2005} for
some relevant results. Recently, Pilehroods proved a $q$-analog of the Two-one formula in \cite{KHP2015}. For instance,
for positive integer $m$, we have
\[\zeta _{}^ \star \left[ {{{\left\{ 2 \right\}}_m}} \right] = \sum\limits_{r = 1}^\infty  {{{\left( { - 1} \right)}^{r - 1}}{q^{r\left( {r + 2m - 1} \right)/2}}} \frac{{\left( {1 + {q^r}} \right)}}
{{\left[ r \right]_q^{2m}}}.\]
Here $\{l\}_m$ denotes the sequence $\underbrace{l,\ldots,l}_{m \text{\;times}}$. By taking $q\rightarrow 1$ in identity above we obtain the following well-known result
\[\zeta _{}^ \star \left( {{{\left\{ 2 \right\}}_m}} \right) = 2\sum\limits_{r = 1}^\infty  {\frac{{{{\left( { - 1} \right)}^{r - 1}}}}
{{{r^{2m}}}}}  = 2\left( {1 - {2^{1 - 2m}}} \right)\zeta \left( {2m} \right).\]

However, we do not consider the $q$-extension of multiple zeta values. We continue the study of $q$-Euler sums in this paper. Our work is motivated by the results recently discovered in \cite{XC2017,X2017}. In \cite{XC2017}, the second author used the integrals
\[\int\limits_0^x {\frac{{{H_m}\left( t,a \right){H_p}\left( t ,b\right)}}{t}} dt, \quad x\in (-1,1)\]
and
\[\int\limits_0^x {\frac{{{H_m}\left( {t,a} \right){H_{p}}\left( {t,b} \right)}}{{t\left( {1 - t} \right)}}dt},\quad x\in (-1,1) \]
to establish many relations involving digamma function, Hurwitz zeta function, parametric linear and quadratic Euler sums. Here the function ${H_k}\left( {x,a} \right)$ is defined by
\[{H_k}\left( {x,a} \right) = \mathop {\lim }\limits_{q \to 1} {H_k}\left[ {x,a} \right].\]
And in \cite{X2017}, we used the Tornheim type series
\[\sum\limits_{n,k = 1}^\infty  {\frac{{H_n^{\left( m \right)}}}{{{k^p}n\left( {n + k} \right)}}},\ k,p,m\in \N\]
to prove the result that the combined quadratic sum
\[{\left( { - 1} \right)^{p - 1}}S\left( {1,m;p + 1} \right) + {\left( { - 1} \right)^{m - 1}}S\left( {1,p + 1;m} \right)\]
are reducible to linear sums and zeta values. We find that the methods and results in \cite{XC2017,X2017} are easily extended to our $q$-Euler sums.
In this paper we will give some extended results on $q$-analogues of Euler sums, see Theorem \ref{Thm:q-Hurwitz}-\ref{Thm:q-Euler-Second}. Moreover, we use certain stuffle products to obtain a general formula of product of any $q$-polylogarithms.

The purpose of the paper is to prove the following four theorems.

\begin{thm}\label{Thm:q-Hurwitz}
Let $k,l$ be positive integers and $a,b,x$ be real numbers with $a,b,a+b\neq -1,-2,\cdots$ and $|x|<1$. Then the following identity holds:
\begin{align}
&(-1)^{k-1}\sum\limits_{m=1}^\infty\frac{q^{(m+b)k}}{[m+b]^{k+l}}\sum\limits_{j=1}^m\frac{x^{j+a+b}}{[j+a+b]}
-(-1)^{l-1}\sum\limits_{m=1}^\infty\frac{q^{(m+a)l}}{[m+a]^{k+l}}\sum\limits_{j=1}^m\frac{x^{j+a+b}}{[j+a+b]}\nonumber\\
& =\sum\limits_{j=1}^{l-1}(-1)^{j-1}H_{k+j}[q^{j-1}x,a]H_{l+1-j}[x,b]
-\sum\limits_{j=1}^{k-1}(-1)^{j-1}H_{l+j}[q^{j-1}x,b]H_{k+1-j}[x,a]\nonumber\\
&\quad +(-1)^{l-1}\left(H_1[x,b]H_{k+l}[q^{l-1}x,a]-H_1[x,a+b]H_{k+l}[q^l,a]\right)\nonumber\\
&\quad -(-1)^{k-1}\left(H_1[x,a]H_{k+l}[q^{k-1}x,b]-H_1[x,a+b]H_{k+l}[q^k,b]\right).
\label{Eq:q-Hurwitz}
\end{align}
\end{thm}

\begin{thm}\label{Thm:q-Euler-First}
Let $k,l$ be positive integers and $s,h,x$ be real numbers with $l>s\geq 0, k>h\geq 0$ and $|x|<1$. Then we have
\begin{align}
&(-1)^{k-1}S\left[ {\left. {\begin{array}{*{20}{c}}
   l,  \\
   q^s, \\
\end{array}\begin{array}{*{20}{c}}
   1  \\
   q^hx  \\
\end{array}} \right|\begin{array}{*{20}{c}}
   k  \\
   q^{k-h}  \\
\end{array}} \right] - (-1)^{l-1}S\left[ {\left. {\begin{array}{*{20}{c}}
   k,  \\
   q^h,  \\
\end{array}\begin{array}{*{20}{c}}
   1  \\
   q^sx  \\
\end{array}} \right|\begin{array}{*{20}{c}}
   l  \\
   q^{l- s}  \\
\end{array}} \right]\nonumber\\
= &\sum\limits_{j=1}^{l-1} (-1)^{j-1}\Li_{l+1-j}[q^sx]S\left[ {\left. {\begin{array}{*{20}{c}}
   k  \\
   {{q^h}}  \\
\end{array}} \right|\begin{array}{*{20}{c}}
   j  \\
   {{q^{j - 1}}x}  \\
\end{array}} \right] - \sum\limits_{j=1}^{k-1} (-1)^{j-1}\Li_{k+1-j}[q^hx]S\left[ {\left. {\begin{array}{*{20}{c}}
   l  \\
   {{q^s}}  \\
\end{array}} \right|\begin{array}{*{20}{c}}
   j  \\
   {{q^{j-1}}x}  \\
\end{array}} \right] \nonumber\\
&\quad + (-1)^{l-1}\ln [1-q^sx]\left( S\left[ \left. {\begin{array}{*{20}{c}}
   k  \\
   q^h  \\
\end{array}} \right|\begin{array}{*{20}{c}}
   l  \\
   q^{l-s}  \\
\end{array} \right] - S\left[ \left. {\begin{array}{*{20}{c}}
   k  \\
   q^h  \\
\end{array}} \right|\begin{array}{*{20}{c}}
   l  \\
   q^{l-1}x  \\
\end{array}\right] \right)\nonumber\\
&\quad - (-1)^{k-1}\ln[1-q^hx]\left(S\left[ \left. {\begin{array}{*{20}{c}}
   l  \\
   q^s  \\
\end{array}} \right|\begin{array}{*{20}{c}}
   k  \\
   q^{k- h}  \\
\end{array}\right] - S\left[ \left. {\begin{array}{*{20}{c}}
   l  \\
   q^s  \\
\end{array}} \right|\begin{array}{*{20}{c}}
   k  \\
   q^{k-1}x  \\
\end{array} \right] \right).
\label{Eq:q-Euler-First}
\end{align}
\end{thm}

\begin{thm}\label{Thm:q-Euler-Second}
For positive integers $k$ and $l$, it holds
\begin{align}
&(-1)^{k-1}S\left[\left.\begin{array}{*{20}{c}}
   l+1,  \\
   q^l,  \\
\end{array}\begin{array}{*{20}{c}}
   1  \\
   q  \\
\end{array} \right|\begin{array}{*{20}{c}}
   k  \\
   q  \\
\end{array} \right]+(-1)^{l-1}S\left[\left.\begin{array}{*{20}{c}}
   k,  \\
   q^{k-1},  \\
\end{array}\begin{array}{*{20}{c}}
   1  \\
   q  \\
\end{array} \right|\begin{array}{*{20}{c}}
   l+1  \\
   q  \\
\end{array} \right]\nonumber\\
=&\Li_{l+1}[q^l]\Li_{k+1}[q^k]+\sum\limits_{j=1}^{k-1}(-1)^{j-1}\Li_{k+1-j}[q^{k-j}]S\left[\left.\begin{array}{*{20}{c}}
   j  \\
   q  \\
\end{array} \right|\begin{array}{*{20}{c}}
   l+1  \\
   q^l  \\
\end{array} \right]\nonumber \\
&\quad +(-1)^{k-1}\Li_{l+1}[q^l]S\left[\left.\begin{array}{*{20}{c}}
   1  \\
   q  \\
\end{array} \right|\begin{array}{*{20}{c}}
   k  \\
   q  \\
\end{array} \right]-\sum\limits_{j=1}^{k-1}(-1)^{j-1}\Li_{k+1-j}[q^{k-j}]\Li_{l+j+1}[q^{l+1}]\nonumber\\
&\quad-\sum\limits_{j=1}^{l-1}(-1)^{j-1}\Li_{l+1-j}[q^{l-j}]S\left[\left.\begin{array}{*{20}{c}}
   k  \\
   q^{k-1}  \\
\end{array} \right|\begin{array}{*{20}{c}}
  j+1 \\
   q  \\
\end{array} \right].
\label{Eq:q-Euler-Second}
\end{align}
\end{thm}

\begin{thm}\label{Thm:q-Polylog-Product}
Let $n,k_1,\ldots,k_n$ be positive integers and $x_1,\ldots,x_n$ be real numbers with $|x_j|<1$. we have
$$\prod\limits_{j=1}^n\Li_{k_j}[x_j]=\sum\limits_{j=0}^{n-1}\sum\limits_{1\leq i_1<\cdots<i_j\leq n}(-1)^{n-1-j}S\left[\left.\begin{array}{*{20}{c}}
   k_{i_1},  \\
   x_{i_1},  \\
\end{array}\begin{array}{*{20}{c}}
   \ldots,  \\
   \ldots,  \\
\end{array}\begin{array}{*{20}{c}}
   k_{i_j}  \\
   x_{i_j}  \\
\end{array} \right|\begin{array}{*{20}{c}}
   (k_1+\cdots+k_n)-(k_{i_1}+\cdots+k_{i_j})  \\
   (x_1\cdots x_n)/(x_{i_1}\cdots x_{i_j})  \\
\end{array} \right]. $$
\end{thm}

We prove Theorems \ref{Thm:q-Hurwitz}-\ref{Thm:q-Euler-Second} in Section \ref{Sec:Proof-1-3} by calculating the Jackson $q$-integral of $q$-polylogarithm functions, and prove Theorem \ref{Thm:q-Polylog-Product} in Section \ref{Sec:Proof-4} algebraically. In Section \ref{Sec:Id-EulerSum}, we give some interesting identities (known or new) involving harmonic numbers.

\section{Proofs of Theorems \ref{Thm:q-Hurwitz}, \ref{Thm:q-Euler-First} and \ref{Thm:q-Euler-Second}}\label{Sec:Proof-1-3}

We prove Theorems \ref{Thm:q-Hurwitz}, \ref{Thm:q-Euler-First} and \ref{Thm:q-Euler-Second} in this section by calculating the Jackson $q$-integral of $q$-polylogarithm functions.

\subsection{Jackson $q$-integral}

The Jackson $q$-integral and $q$-derivative are defined by (\cite{A2000,AM2012,Ba2014,J2005,Xu2016})
\begin{align*}
&\int\limits_a^x {f(t){d_q}t  } := (1-q)\sum\limits_{i = 0}^\infty  {{q^i}\left[ {xf(q^ix) - af(q^ia)} \right]},\\
&{D_q}f(x) := \frac{f(qx)-f(x)}{qx - x},
\end{align*}
respectively. For example, we have
$$D_q(H_k[x,a]) = \frac{H_{k-1}[x,a]}{x},\qquad {D_q}(x^m) = [m]x^{m- 1}.$$
And it is easy to verify that
\begin{align*}
&D_q(f(x)g(x))=g(qx)D_q(f(x))+f(x)D_q(g(x))= f(qx)D_q(g(x))+g(x)D_q(f(x)),\\
&D_q\left(\int\limits_a^x f(t)d_qt\right)=f(x),\qquad\int\limits_a^x D_q(f(t))d_qt=f(x)-f(a),\\
&\int\limits_a^x f(t)D_q(g(t))d_qt=\left.\left[f(t)g(t)\right]\right|_a^x-\int\limits_a^x g(qt)D_q(f(t))d_qt.
\end{align*}

\subsection{Proof of Theorem \ref{Thm:q-Hurwitz}}

To prove Theorem \ref{Thm:q-Hurwitz}, we need a lemma.

\begin{lem}
Let $m,k$ be positive integers and $a,b,x$ be real numbers with $a,b,a+b\neq -1,-2,\cdots$ and $|x|<1$. Then the following identity holds:
\begin{align}
\int\limits_0^x H_k[t,a]t^{m+b-1}d_qt=&\sum\limits_{j=1}^{k-1}(-1)^{j-1}\frac{q^{(m+b)(j-1)}x^{m+b}}{[m+b]^j}H_{k+1-j}[x,a]\nonumber\\
&+(-1)^{k-1}\frac{q^{(m+b)(k-1)}}{[m+b]^k}\left(x^{m+b}H_1[x,a]-q^{m+b}H_1[x,a+b]\right)\nonumber\\
&+(-1)^{k-1}\frac{q^{(m+b)k}}{[m+b]^k}\sum\limits_{j=1}^m\frac{x^{j+a+b}}{[j+a+b]}.
\label{Eq:Jackson-Hirwitz}
\end{align}
\end{lem}

\pf
Denote the left hand-side of \eqref{Eq:Jackson-Hirwitz} by $I_k$. Then we have
\begin{align*}
I_k=\frac{1}{[m+b]}\int\limits_0^x H_k[t,a]D_q(t^{m+b})d_qt=\frac{x^{m+b}}{[m+b]}H_k[x,a]-\frac{q^{m+b}}{[m+b]}I_{k-1},
\end{align*}
and
\begin{align*}
I_1=&\frac{x^{m+b}}{[m+b]}H_1[x,a]-\frac{q^{m+b}}{[m+b]}\int\limits_{0}^x\frac{t^{m+a+b}}{1-t}d_qt\\
=&\frac{x^{m+b}}{[m+b]}H_1[x,a]+\frac{q^{m+b}}{[m+b]}\sum\limits_{j=1}^m\frac{x^{j+a+b}}{[j+a+b]}-\frac{q^{m+b}}{[m+b]}H_1[x,a+b].
\end{align*}
Hence we get \eqref{Eq:Jackson-Hirwitz} by induction on $k$.
\qed

\noindent {\bf Proof of Theorem \ref{Thm:q-Hurwitz}.}
Considering the Jackson $q$-integral
\begin{align*}
\int\limits_0^x \frac{H_k[t,a]H_l[t,b]}{t}d_qt =&\sum\limits_{m=1}^\infty\frac{1}{[m+a]^k}\int\limits_0^x H_l[t,b]t^{m+a-1}d_qt \\
=&\sum\limits_{m=1}^\infty  \frac{1}{[m+b]^l}\int\limits_0^x H_k[t,a]t^{m+b-1}d_qt,
\end{align*}
we get \eqref{Eq:q-Hurwitz} with the help of \eqref{Eq:Jackson-Hirwitz}.
\qed

Setting $x=q$ and $k=l=1$ in Theorem \ref{Thm:q-Hurwitz}, we get
\begin{align*}
&\sum\limits_{m=1}^\infty  \frac{q^{m+b}}{[m+b]^2}\sum\limits_{j=1}^m\frac{q^{j+a+b}}{[j+a+b]}
-\sum\limits_{m=1}^\infty  \frac{q^{m+a}}{[m+a]^2}\sum\limits_{j=1}^m \frac{q^{j+a+b}}{[j+a+b]}\\
=& [a]H_2[q,a]\sum\limits_{m=1}^\infty \frac{q^{m+b}}{[m+b][m+a+b]}
-[b]H_2[q,b]\sum\limits_{m=1}^\infty \frac{q^{m+a}}{[m+a][m+a+b]}.
\end{align*}

\subsection{Proof of Theorem \ref{Thm:q-Euler-First}}

Similarly as the proof of Theorem \ref{Thm:q-Hurwitz}, we give a proof of Theorem \ref{Thm:q-Euler-First}.

\noindent {\bf Proof of Theorem \ref{Thm:q-Euler-First}.}
Using \eqref{Eq:Jackson-Hirwitz}, we compute the Jackson $q$-integral
\begin{align*}
&\int\limits_0^x \frac{\Li_l[q^st]\Li_{k}[q^ht]}{t(1-t)}d_qt
=\sum\limits_{m=1}^\infty  \zeta_m[l,q^s]\int\limits_0^x t^{m-1}\Li_k[q^ht]d_qt \\
=&\sum\limits_{j=1}^{k-1}(-1)^{j-1}\Li_{k+1-j}[q^hx]\sum\limits_{m=1}^\infty\frac{\zeta_m[l,q^s]}{[m]^j}\left(q^{j-1}x\right)^m\\
&+(-1)^{k-1}\ln[1-q^hx]\sum\limits_{m=1}^\infty\frac{\zeta_m[l,q^s]}{[m]^k}\left(q^{(k-h)m}-(q^{k-1}x)^m\right)\\
&+(-1)^{k-1}\sum\limits_{m=1}^\infty\frac{\zeta_m[l,q^s]\zeta_m[1,q^hx]}{[m]^k}q^{(k-h)m}\\
=&\sum\limits_{j=1}^{l-1}(-1)^{j-1}\Li_{l+1-j}[q^sx]\sum\limits_{m=1}^\infty\frac{\zeta_m[k,q^h]}{[m]^j}\left(q^{j-1}x\right)^m\\
&+(-1)^{l-1}\ln[1-q^sx]\sum\limits_{m=1}^\infty\frac{\zeta_m[k,q^h]}{[m]^l}\left(q^{(l-s)m}-(q^{l-1}x)^m\right)\\
&+(-1)^{l-1}\sum\limits_{m=1}^\infty\frac{\zeta_m[k,q^h]\zeta_m[1,q^sx]}{[m]^l}q^{(l-s)m},
\end{align*}
from which we get \eqref{Eq:q-Euler-First}.\qed

Setting $x\rightarrow 1$ in Theorem \ref{Thm:q-Euler-First}, we obtain
\begin{align}
&(-1)^{k-1}S\left[ \left. {\begin{array}{*{20}{c}}
   l,  \\
   q^s,  \\
\end{array}\begin{array}{*{20}{c}}
   1  \\
   q^h  \\
\end{array}} \right|\begin{array}{*{20}{c}}
   k  \\
   q^{k - h}  \\
\end{array} \right] - (-1)^{l-1}S\left[ \left. {\begin{array}{*{20}{c}}
   k,  \\
   q^h,  \\
\end{array}\begin{array}{*{20}{c}}
   1  \\
   q^s  \\
\end{array}} \right|\begin{array}{*{20}{c}}
   l  \\
   q^{l-s}  \\
\end{array} \right]\nonumber\\
=&\sum\limits_{j=2}^{l-1} (-1)^{j-1}\Li_{l+1-j}[q^s]S\left[ \left. {\begin{array}{*{20}{c}}
   k  \\
   q^h  \\
\end{array}} \right|\begin{array}{*{20}{c}}
   j  \\
   q^{j-1}  \\
\end{array} \right]-\sum\limits_{j=2}^{k-1} (-1)^{j-1}\Li_{k+1-j}[q^h]S\left[ \left. {\begin{array}{*{20}{c}}
   l  \\
   q^s  \\
\end{array}} \right|\begin{array}{*{20}{c}}
   j  \\
   q^{j-1}  \\
\end{array} \right] \nonumber\\
&\quad +(-1)^{l-1}\ln[1-q^s]\left(S\left[ \left. \begin{array}{*{20}{c}}
   k  \\
   q^h  \\
\end{array} \right|\begin{array}{*{20}{c}}
   l  \\
   q^{l-s}  \\
\end{array} \right] - S\left[ \left. \begin{array}{*{20}{c}}
   k  \\
   q^h  \\
\end{array} \right|\begin{array}{*{20}{c}}
   l  \\
   q^{l-1}  \\
\end{array} \right] \right)\nonumber\\
&\quad-(-1)^{k-1}\ln[1-q^h]\left(S\left[\left. \begin{array}{*{20}{c}}
   l  \\
   q^s  \\
\end{array} \right|\begin{array}{*{20}{c}}
   k  \\
   q^{k- h}  \\
\end{array} \right] - S\left[ \left. \begin{array}{*{20}{c}}
   l  \\
   q^s  \\
\end{array} \right|\begin{array}{*{20}{c}}
   k  \\
   q^{k-1}  \\
\end{array} \right] \right)\nonumber\\
&\quad +\sum\limits_{m=1}^\infty\frac{\Li_l[q^s]\zeta_m[k,q^h]-\Li_k[q^h]\zeta_m[l,q^s]}{[m]}.
\label{Eq:q-Euler-First-x1}
\end{align}

To evaluate the last sum in the right-hand side of \eqref{Eq:q-Euler-First-x1}, we use

\begin{thm}\label{Thm:Mul-Li-Zeta}
Let $k,l$ be positive integers and $x,y,z$ be real numbers with $|x|,|y|,|z|<1$. Then we have
\begin{align}
&\sum\limits_{m=1}^\infty \frac{\Li_k[x]\zeta_m[l,y]-\Li_l[y]\zeta_m[k,x]}{[m]}z^m  \nonumber\\
=&\Li_l[y]S\left[\left.\begin{array}{*{20}{c}}
   1  \\
   z  \\
\end{array} \right|\begin{array}{*{20}{c}}
   k  \\
   x  \\
\end{array} \right]-\Li_k[x]S\left[\left.\begin{array}{*{20}{c}}
   1  \\
   z  \\
\end{array} \right|\begin{array}{*{20}{c}}
   l  \\
   y  \\
\end{array} \right]+\Li_k[x]\Li_{l+1}[zy]-\Li_l[y]\Li_{k+1}[zx].
\label{Eq:Mul-Li-Zeta}
\end{align}
\end{thm}

Taking $x=q^h,y=q^s$ and $z\rightarrow1$ in \eqref{Eq:Mul-Li-Zeta}, and using \eqref{Eq:q-Euler-First-x1}, we have

\begin{cor}
Let $s,h$ be positive reals and $k,l$ be positive integers with $l>\max\{s,1\}$ and $k>\max\{h,1\}$. Then it holds
\begin{align*}
&(-1)^{k-1}S\left[ \left. {\begin{array}{*{20}{c}}
   l,  \\
   q^s,  \\
\end{array}\begin{array}{*{20}{c}}
   1  \\
   q^h  \\
\end{array}} \right|\begin{array}{*{20}{c}}
   k  \\
   q^{k - h}  \\
\end{array} \right] - (-1)^{l-1}S\left[ \left. {\begin{array}{*{20}{c}}
   k,  \\
   q^h,  \\
\end{array}\begin{array}{*{20}{c}}
   1  \\
   q^s  \\
\end{array}} \right|\begin{array}{*{20}{c}}
   l  \\
   q^{l-s}  \\
\end{array} \right]\\
=&\sum\limits_{j=2}^{l-1} (-1)^{j-1}\Li_{l+1-j}[q^s]S\left[ \left. {\begin{array}{*{20}{c}}
   k  \\
   q^h  \\
\end{array}} \right|\begin{array}{*{20}{c}}
   j  \\
   q^{j-1}  \\
\end{array} \right]-\sum\limits_{j=2}^{k-1} (-1)^{j-1}\Li_{k+1-j}[q^h]S\left[ \left. {\begin{array}{*{20}{c}}
   l  \\
   q^s  \\
\end{array}} \right|\begin{array}{*{20}{c}}
   j  \\
   q^{j-1}  \\
\end{array} \right] \\
&\quad +(-1)^{l-1}\ln[1-q^s]\left(S\left[ \left. \begin{array}{*{20}{c}}
   k  \\
   q^h  \\
\end{array} \right|\begin{array}{*{20}{c}}
   l  \\
   q^{l-s}  \\
\end{array} \right] - S\left[ \left. \begin{array}{*{20}{c}}
   k  \\
   q^h  \\
\end{array} \right|\begin{array}{*{20}{c}}
   l  \\
   q^{l-1}  \\
\end{array} \right] \right)\\
&\quad-(-1)^{k-1}\ln[1-q^h]\left(S\left[\left. \begin{array}{*{20}{c}}
   l  \\
   q^s  \\
\end{array} \right|\begin{array}{*{20}{c}}
   k  \\
   q^{k- h}  \\
\end{array} \right] - S\left[ \left. \begin{array}{*{20}{c}}
   l  \\
   q^s  \\
\end{array} \right|\begin{array}{*{20}{c}}
   k  \\
   q^{k-1}  \\
\end{array} \right] \right)\\
&\quad +\Li_k[q^h]S\left[\left.\begin{array}{*{20}{c}}
   1  \\
   1  \\
\end{array} \right|\begin{array}{*{20}{c}}
   l  \\
   q^s  \\
\end{array} \right]-\Li_l[q^s]S\left[\left.\begin{array}{*{20}{c}}
   1  \\
   1  \\
\end{array} \right|\begin{array}{*{20}{c}}
   k  \\
   q^h  \\
\end{array} \right]+\Li_l[q^s]\Li_{k+1}[q^h]-\Li_k[q^h]\Li_{l+1}[q^s].
\end{align*}
\end{cor}

Finally, we give a proof of Theorem \ref{Thm:Mul-Li-Zeta}.

\noindent {\bf Proof of Theorem \ref{Thm:Mul-Li-Zeta}.}
We compute the $N$-th partial sum of the series of the left-hand side of \eqref{Eq:Mul-Li-Zeta}
\begin{align*}
&\sum\limits_{m=1}^N \frac{\Li_k[x]\zeta_m[l,y]-\Li_l[y]\zeta_m[k,x]}{[m]}z^m  \\
=&\Li_k[x]\sum\limits_{m=1}^N \frac{\zeta_m[l,y]}{[m]}z^m-\Li_l[y]\sum\limits_{m=1}^N\frac{\zeta_m[k,x]}{[m]}z^m  \\
=&\Li_k[x]\sum\limits_{m=1}^N\sum\limits_{j=1}^m\frac{y^jz^m}{[m][j]^l}-\Li_l[y]\sum\limits_{m=1}^N\sum\limits_{j=1}^m\frac{x^jz^m}{[m][j]^k}  \\
=&\Li_k[x]\sum\limits_{j=1}^N\sum\limits_{m=j}^N\frac{y^jz^m}{[m][j]^l}-\Li_l[y]\sum\limits_{j=1}^N\sum\limits_{m=j}^N\frac{x^jz^m}{[m][j]^k}\\
=&\Li_k[x]\sum\limits_{j=1}^N\frac{\zeta_N[1,z]-\zeta_{j-1}[1,z]}{[j]^l}y^j-\Li_l[y]\sum\limits_{j=1}^N\frac{\zeta_N[1,z]-\zeta_{j-1}[1,z]}{[j]^k}x^j\\
=&\zeta_N[1,z](\Li_k[x]\zeta_N[l,y]-\Li_l[y]\zeta_N[k,x])+\Li_l[y]\sum\limits_{j=1}^N\frac{\zeta_{j-1}[1,z]}{[j]^k}x^j-\Li_k[x]\sum\limits_{j=1}^N\frac{\zeta_{j-1}[1,z]}{[j]^l}y^j.
\end{align*}
Letting $N$ tend to infinity, we get \eqref{Eq:Mul-Li-Zeta}.
\qed

\subsection{Proof of Theorem \ref{Thm:q-Euler-Second}}

To prove Theorem \ref{Thm:q-Euler-Second}, we need the following lemmas.

\begin{lem}
For positive integers $k$ and $i$, it holds
\begin{align}
\sum\limits_{m=1}^\infty \frac{\zeta_m[k,q^{k-1}]}{[m][m+i]}q^m=\frac{1}{[i]}\left\{\begin{array}{l}
 \Li_{k+1}[q^k]+(-1)^{k-1}\sum\limits_{j=1}^{i-1}\frac{\left[H_j^{(1)}\right]}{[j]^k}q^j \\
  +\sum\limits_{j=1}^{k-1}(-1)^{j-1}\Li_{k+1-j}[q^{k-j}]\left[H_{i-1}^{(j)}\right]
 \end{array} \right\},
\label{Eq:sum-Zeta}
\end{align}
where we set \cite{Xu2016}
$$\left[H_m^{(k)}\right]=\zeta_m[k,q].$$
\end{lem}

\pf By using the Cauchy product of power series and the definition of $q$-harmonic numbers, we have
$$\sum\limits_{m=1}^\infty\zeta_m[k,q^l]x^m=\frac{\Li_k[q^lx]}{1-x},$$
where $l$ is any positive integer and $x$ is any real number with $|x|<1$.
Multiplied by $x^{-1}-x^{i-1}$ and $q$-integrated over $(0,q)$, the above equation yields
$$[i]\sum\limits_{m=1}^\infty\frac{\zeta_m[k,q^l]}{[m][m+i]}q^m=\Li_{k+1}[q^{l+1}]+\sum\limits_{j=1}^{i-1}\sum\limits_{m=1}^\infty\frac{q^{(l+1)m+j}}{[m]^k[m+j]}.$$
Taking $l=k-1$ in above equation, we obtain
$$[i]\sum\limits_{m=1}^\infty\frac{\zeta_m[k,q^{k-1}]}{[m][m+i]}q^m=\Li_{k+1}[q^{k}]+\sum\limits_{j=1}^{i-1}q^j\sum\limits_{m=1}^\infty\frac{q^{km}}{[m]^k[m+j]},$$
which together with the formula
\begin{align}
\sum\limits_{m=1}^\infty\frac{q^{km}}{[m]^k[m+j]}=\sum\limits_{p=1}^{k-1}\frac{(-1)^{p-1}}{[j]^p}\Li_{k-p+1}[q^{k-p}]+(-1)^{k-1}\frac{\left[H_j^{(1)}\right]}{[j]^k}
\label{Eq:Double-Sum}
\end{align}
yield the desired result. \qed

\begin{lem}
For any positive integers $k_1,k_2$ and any real numbers $x,y$ with $|x|,|y|<1$, we have
\begin{align}
S\left[\left.\begin{array}{*{20}{c}}
   k_1  \\
   x  \\
\end{array} \right|\begin{array}{*{20}{c}}
   k_2  \\
   y  \\
\end{array} \right] + S\left[\left.\begin{array}{*{20}{c}}
   k_2 \\
   y  \\
\end{array} \right|\begin{array}{*{20}{c}}
   k_1  \\
   x  \\
\end{array} \right]=\Li_{k_1}[x]\Li_{k_2}[y]+\Li_{k_1+k_2}[xy].
\label{Eq:Li-Li-2}
\end{align}
\end{lem}

\pf We consider the generating function
$$F_2[x,y,z]:=\sum\limits_{m=1}^\infty \left(\zeta_m[k_1,x]\zeta_m[k_2,y]-\zeta_m[k_1+k_2,xy]\right)z^{m-1},$$
where $|z|<1$.
By the definition of $\zeta_m[k,x]$, we have
\begin{align*}
F_2[x,y,z]&=\sum\limits_{m=1}^\infty\left\{\begin{array}{l}
 \left(\zeta_m[k_1,x]+\frac{x^{m+1}}{[m+1]^{k_1}} \right)\left(\zeta_m[k_2,y]+\frac{y^{m+1}}{[m+1]^{k_2}}\right) \\
  -\left(\zeta_m[k_1+k_2,xy]+\frac{x^{m+1}y^{m+1}}{[m+1]^{k_1+k_2}}\right)
 \end{array} \right\}z^{m} \\
=&zF_2[x,y,z]+\sum\limits_{m=1}^\infty\left(\frac{\zeta_m[k_1,x]}{[m+1]^{k_2}}y^{m+1}+\frac{\zeta_m[k_2,y]}{[m+1]^{k_1}}x^{m+1}\right)z^m\\
=&zF_2[x,y,z]+\sum\limits_{m=1}^\infty\left(\frac{\zeta_m[k_1,x]}{[m]^{k_2}}y^m+\frac{\zeta_m[k_2,y]}{[m]^{k_1}}x^m-2\frac{x^my^m}{[m]^{k_1+k_2}}\right)z^{m-1}.
 \end{align*}
Hence, we obtain
\begin{align*}
F_2[x,y,z]=&\sum\limits_{m=1}^\infty\left(\frac{\zeta_m[k_1,x]}{[m]^{k_2}}y^m+\frac{\zeta_m[k_2,y]}{[m]^{k_1}}x^m-2\frac{x^my^m}{[m]^{k_1+k_2}}\right)\frac{z^{m-1}}{1-z}\\
=& \sum\limits_{m=1}^\infty\sum\limits_{j=1}^m\left(\frac{\zeta_j[k_1,x]}{[j]^{k_2}}y^j+\frac{\zeta_j[k_2,y]}{[j]^{k_1}}x^j-2\frac{x^jy^j}{[j]^{k_1+k_2}}\right)z^{m-1}.
\end{align*}
Then equating coefficients of $z^{m-1}$, we establish the relation
$$\sum\limits_{j=1}^m\left(\frac{\zeta_j[k_1,x]}{[j]^{k_2}}y^j+\frac{\zeta_j[k_2,y]}{[j]^{k_1}}x^j\right)=\zeta_m[k_1,x]\zeta_m[k_2,y]+\zeta_m[k_1+k_2,xy].$$
Letting $m$ tend to infinity in above equation, we deduce \eqref{Eq:Li-Li-2}.\qed

\begin{re}
Similarly, considering the following function
$$F_3[x,y,z,t]:=\sum\limits_{m=1}^\infty\left(\zeta_m[k_1,x]\zeta_m[k_2,y]\zeta_m[k_3,z]-\zeta_m[k_1+k_2+k_3,xyz]\right)t^{m-1},$$
and applying the same arguments as in the proof of \eqref{Eq:Li-Li-2}, we may deduce the following formula
\begin{align}
&S\left[\left.\begin{array}{*{20}{c}}
   k_1,  \\
   x,  \\
\end{array}\begin{array}{*{20}{c}}
   k_2  \\
   y  \\
\end{array} \right|\begin{array}{*{20}{c}}
   k_3  \\
   z  \\
\end{array} \right] + S\left[\left.\begin{array}{*{20}{c}}
   k_1,  \\
   x,  \\
\end{array}\begin{array}{*{20}{c}}
   k_3  \\
   z  \\
\end{array} \right|\begin{array}{*{20}{c}}
   k_2  \\
   y  \\
\end{array} \right] + S\left[\left.\begin{array}{*{20}{c}}
   k_2,  \\
   y,  \\
\end{array}\begin{array}{*{20}{c}}
   k_3  \\
   z  \\
\end{array} \right|\begin{array}{*{20}{c}}
   k_1  \\
   x  \\
\end{array} \right]\nonumber\\
=& S\left[ \left. \begin{array}{*{20}{c}}
   k_1  \\
   x  \\
\end{array} \right|\begin{array}{*{20}{c}}
   k_2+k_3  \\
   yz  \\
\end{array} \right] + S\left[\left.\begin{array}{*{20}{c}}
   k_2  \\
   y  \\
\end{array} \right|\begin{array}{*{20}{c}}
   k_1+k_3  \\
   xz  \\
\end{array} \right] + S\left[\left.\begin{array}{*{20}{c}}
   k_3  \\
   z  \\
\end{array} \right|\begin{array}{*{20}{c}}
   k_1+k_2  \\
   xy  \\
\end{array} \right]\nonumber\\
&\quad+\Li_{k_1}[x]\Li_{k_2}[y]\Li_{k_3}[z]-\Li_{k_1+k_2+k_3}[xyz],
\label{Eq:Li-Li-3}
\end{align}
where $k_1,k_2,k_3$ are positive integers and $x,y,z$ are real number with $|x|,|y|,|z|<1$. It is clear that \eqref{Eq:Li-Li-2} and \eqref{Eq:Li-Li-3} are immediate corollaries of Theorem \ref{Thm:q-Polylog-Product}.
\end{re}

\begin{lem}
For any positive integers $k_1,k_2$ and any real numbers $x,y,z$ with $|x|,|y|,|z|<1$, we have
\begin{align}
\sum\limits_{m=1}^\infty\frac{x^m}{[m]^{k_1}}\sum\limits_{j=1}^m\frac{z^j}{[j]^{k_2}}\zeta_j[1,y]+S\left[\left.\begin{array}{*{20}{c}}
   k_1,  \\
   x,  \\
\end{array}\begin{array}{*{20}{c}}
   1  \\
   y  \\
\end{array} \right|\begin{array}{*{20}{c}}
   k_2  \\
   z  \\
\end{array} \right] = \Li_{k_1}[x]S\left[\left.\begin{array}{*{20}{c}}
   1  \\
   y  \\
\end{array} \right|\begin{array}{*{20}{c}}
   k_2  \\
   z  \\
\end{array} \right] + S\left[\left.\begin{array}{*{20}{c}}
   1  \\
   y  \\
\end{array} \right|\begin{array}{*{20}{c}}
   k_1+k_2  \\
   xz  \\
\end{array}\right].
\label{Eq:Li-S}
\end{align}
\end{lem}

\pf Replacing $y$ by $zt$ in \eqref{Eq:Li-Li-2}, then dividing it by $\frac{1}{1-t}$ and $q$-integrating over the interval $(0,y)$, we can deduce the desired result.\qed

Finally, we come to the proof of Theorem \ref{Thm:q-Euler-Second}.

\noindent {\bf Proof of Theorem \ref{Thm:q-Euler-Second}.}
Set
$$\sum=\sum\limits_{i=1}^\infty\sum\limits_{m=1}^\infty\frac{\zeta_m[k,q^{k-1}]}{[m][m+i][i]^l}q^{m+li}.$$
On the one hand, using \eqref{Eq:sum-Zeta}, we have
\begin{align*}
\sum=&\sum\limits_{i=1}^\infty\frac{q^{li}}{[i]^{l+1}}\left\{\begin{array}{l}
 \Li_{k+1}[q^k]+(-1)^{k-1}\sum\limits_{j=1}^{i-1}\frac{\left[H_j^{(1)}\right]}{[j]^k}q^j \\
  +\sum\limits_{j=1}^{k-1}(-1)^{j-1}\Li_{k+1-j}[q^{k-j}]\left[H_{i-1}^{(j)}\right]
 \end{array} \right\} \\
=& (-1)^{k-1}\sum\limits_{m=1}^\infty\frac{q^{lm}}{[m]^{l+1}}\sum\limits_{j=1}^m\frac{\left[H_j^{(1)}\right]}{[j]^k}q^j-(-1)^{k-1}S\left[\left.\begin{array}{*{20}{c}}
   1  \\
   q  \\
\end{array} \right|\begin{array}{*{20}{c}}
   k+l+1  \\
   q^{l+1}  \\
\end{array} \right]\\
&\quad +\sum\limits_{j=1}^{k-1}(-1)^{j-1}\Li_{k+1-j}[q^{k-j}]S\left[\left.\begin{array}{*{20}{c}}
   j  \\
   q  \\
\end{array} \right|\begin{array}{*{20}{c}}
   l+1  \\
   q^l  \\
\end{array} \right]+\Li_{l+1}[q^l]\Li_{k+1}[q^k]\\
&\quad-\sum\limits_{j=1}^{k-1}(-1)^{j-1}\Li_{k+1-j}[q^{k-j}]\Li_{l+j+1}[q^{l+1}].
\end{align*}
Setting $k_1=l+1,k_2=k,x=q^l,y=z=q$ in \eqref{Eq:Li-S}, we get
$$\sum\limits_{m=1}^\infty\frac{q^{lm}}{[m]^{l+1}}\sum\limits_{j=1}^m\frac{q^j}{[j]^k}\left[H_j^{(1)}\right]-S\left[\left.\begin{array}{*{20}{c}}
   1  \\
   q  \\
\end{array} \right|\begin{array}{*{20}{c}}
   k+l+1  \\
   q^{l+1}  \\
\end{array} \right]=\Li_{l+1}[q^l]S\left[\left.\begin{array}{*{20}{c}}
   1  \\
   q  \\
\end{array} \right|\begin{array}{*{20}{c}}
   k  \\
   q  \\
\end{array}\right]-S\left[\left.\begin{array}{*{20}{c}}
   l+1,  \\
   q^l,  \\
\end{array}\begin{array}{*{20}{c}}
   1  \\
   q  \\
\end{array} \right|\begin{array}{*{20}{c}}
   k  \\
   q  \\
\end{array} \right],$$
which deduce that
\begin{align}
\sum=&(-1)^{k-1}\Li_{l+1}[q^l]S\left[\left.\begin{array}{*{20}{c}}
   1  \\
   q  \\
\end{array} \right|\begin{array}{*{20}{c}}
   k  \\
   q  \\
\end{array}\right]-(-1)^{k-1}S\left[\left.\begin{array}{*{20}{c}}
   l+1,  \\
   q^l,  \\
\end{array}\begin{array}{*{20}{c}}
   1  \\
   q  \\
\end{array} \right|\begin{array}{*{20}{c}}
   k  \\
   q  \\
\end{array} \right]\nonumber\\
&\quad +\sum\limits_{j=1}^{k-1}(-1)^{j-1}\Li_{k+1-j}[q^{k-j}]S\left[\left.\begin{array}{*{20}{c}}
   j  \\
   q  \\
\end{array} \right|\begin{array}{*{20}{c}}
   l+1  \\
   q^l  \\
\end{array} \right]+\Li_{l+1}[q^l]\Li_{k+1}[q^k]\nonumber\\
&\quad-\sum\limits_{j=1}^{k-1}(-1)^{j-1}\Li_{k+1-j}[q^{k-j}]\Li_{l+j+1}[q^{l+1}].
\label{Eq:Sum-1}
\end{align}

On the other hand, using \eqref{Eq:Double-Sum}, we get
\begin{align}
\sum=& \sum\limits_{m=1}^\infty\frac{\zeta_m[k,q^{k-1}]}{[m]}q^m\sum\limits_{i=1}^\infty\frac{q^{li}}{[i]^l[m+i]}\nonumber\\
=&\sum\limits_{m=1}^\infty\frac{\zeta_m[k,q^{k-1}]}{[m]}q^m\left(\sum\limits_{j=1}^{l-1}\frac{(-1)^{j-1}}{[m]^j}\Li_{l-j+1}[q^{l-j}]
+(-1)^{l-1}\frac{\left[H_{m}^{(1)}\right]}{[m]^l}\right)\nonumber\\
=&\sum\limits_{j=1}^{l-1}(-1)^{j-1}\Li_{l+1-j}[q^{l-j}]S\left[\left.\begin{array}{*{20}{c}}
   k  \\
   q^{k-1}  \\
\end{array} \right|\begin{array}{*{20}{c}}
   j+1  \\
   q  \\
\end{array} \right]+(-1)^{l-1}S\left[\left.\begin{array}{*{20}{c}}
   k,  \\
   q^{k-1},  \\
\end{array}\begin{array}{*{20}{c}}
   1  \\
   q  \\
\end{array} \right|\begin{array}{*{20}{c}}
   l+1  \\
   q  \\
\end{array} \right].
\label{Eq:Sum-2}
\end{align}
Comparing \eqref{Eq:Sum-1} and \eqref{Eq:Sum-2}, we get the result.\qed

\section{Proof of Theorem \ref{Thm:q-Polylog-Product}}\label{Sec:Proof-4}

In this section, we use the stuffle product to give a proof of Theorem \ref{Thm:q-Polylog-Product}.

For a sequence $\mathbf{k}=(k_1,\ldots,k_n)$ of positive integers, the weight and the depth of $\mathbf{k}$ are defined by
$${\rm wt}(\mathbf{k})=k_1+\cdots+k_n,\qquad {\rm dep}(\mathbf{k})=n,$$
respectively. For an empty sequence, we set ${\rm wt}(\varnothing)={\rm dep}(\varnothing)=0$. We call $\mathbf{l}$ a subsequence of $\mathbf{k}$, if there exist integers $m,i_1,\ldots,i_m$ with $0\leqslant m\leqslant n$ and $1\leqslant i_1<\cdots<i_m\leqslant n$, such that $\mathbf{l}=(k_{i_1},\ldots,k_{i_m})$. Let ${\rm Sub}(\mathbf{k})$ be the set of all subsequences of $\mathbf{k}$. If $\mathbf{x}=(x_1,\ldots,x_n)$ is a sequence of variables, we set $|\mathbf{x}|=x_1\cdots x_n$. And for any $\mathbf{l}=(k_{i_1},\ldots,k_{i_m})\in{\rm Sub}(\mathbf{k})$, we set
$$\mathbf{x}_{\mathbf{l}}=(x_{i_1},\ldots,x_{i_m}).$$
Note that $|\varnothing|=1$, $\mathbf{x}_{\varnothing}=\varnothing$ and $\mathbf{x}_{\mathbf{k}}=\mathbf{x}$.
Therefore Theorem \ref{Thm:q-Polylog-Product} may be rewritten as

\begin{thm}\label{Thm:q-Polylog-Product-New}
Let $n$ be a positive integer, $\mathbf{k}=(k_1,\ldots,k_n)$ be a sequence of positive integers and $\mathbf{x}=(x_1,\ldots,x_n)$ be a sequence of real numbers with $|x_j|<1$. we have
\begin{align}
\prod\limits_{j=1}^n \Li_{k_j}[x_j]=\sum\limits_{{\bf{l}}\in {\rm Sub}({\bf{k}}),{\bf{l}} \ne {\bf{k}}} (-1)^{n-{\rm dep}(\bf{l})-1}S\left[\left. \begin{array}{*{20}{c}}
   \bf{l}  \\
   \bf{x}_{\bf{l}}  \\
\end{array} \right|\begin{array}{*{20}{c}}
   {\rm wt}(\bf{k}) - {\rm wt}(\bf{l})  \\
   |\bf{x}|/|\bf{x}_{\bf{l}}|  \\
\end{array} \right].
\label{Eq:q-Polylog-Product}
\end{align}
\end{thm}

To prove Theorem \ref{Thm:q-Polylog-Product-New}, we use the stuffle product. Similar as in \cite{GZ2010,H1997}, let
$$\mathcal{M}:=\left\{\left.\left[\begin{matrix}
k\\
x
\end{matrix}\right]\right|(k,x)\in\mathbb{N}\times (-1,1)\right\},$$
which we regard as an alphabet with noncommutative letters. Let $\mathcal{M}^{\ast}$ be the set of all words generated by $\mathcal{M}$, which contains the empty word $1_{\mathcal{M}}$. We denote a nonempty word $\left[\begin{matrix}
k_1\\
x_1
\end{matrix}\right]\cdots\left[\begin{matrix}
k_n\\
x_n
\end{matrix}\right]$ simplify by $\left[\begin{matrix}
k_1,\ldots,k_n\\
x_1,\ldots,x_n
\end{matrix}\right]$. Let $\mathfrak{h}^1=\mathbb{Q}\langle\mathcal{M}\rangle$ be the noncommutative polynomial algebra over $\mathbb{Q}$ generated by $\mathcal{M}$. As a rational vector space, $\mathfrak{h}^1$ has a basis $\mathcal{M}^{\ast}$.

We now define the stuffle product ${\bar \ast}$ on the algebra $\mathfrak{h}^1$, which is $\mathbb{Q}$-bilinear, and satisfies the following axioms
\begin{itemize}
  \item [(1)] $1_\mathcal{M}{\bar \ast} w=w{\bar \ast} 1_\mathcal{M}=w$ for any $w\in\mathcal{M}^{\ast}$;
  \item [(2)] $au{\bar \ast} bv=a(u{\bar \ast} bv)+b(au{\bar \ast} v)-(a\circ b)(u{\bar \ast} v)$ for any $a,b\in\mathcal{M}$ and any $u,v\in\mathcal{M}^{\ast}$.
\end{itemize}
Here we set
$$\left[\begin{matrix}
k\\
x
\end{matrix}\right]\circ\left[\begin{matrix}
l\\
y
\end{matrix}\right]:=\left[\begin{matrix}
k+l\\
xy
\end{matrix}\right].$$
Then by \cite{H1997,M2009}, the product ${\bar \ast}$ is commutative and associative.

For any $w\in\mathfrak{h}^1$, we define a function ${\rm Li}^{\star}[w]$ by $\mathbb{Q}$-linearity, ${\rm Li}^{\star}[1_{\mathcal{M}}]=1$ and
$${\rm Li}^{\star}\left[\begin{matrix}
k_1, \ldots,k_n\\
x_1,\ldots,x_n
\end{matrix}\right]:=\sum\limits_{m_1\geqslant \cdots\geqslant m_n\geqslant 1}\frac{x_1^{m_1}\cdots x_n^{m_n}}{[m_1]^{k_1}\cdots[m_n]^{k_n}}.$$
Then we have
$${\rm Li}^{\star}\left[\begin{matrix}
k\\
x
\end{matrix}\right]={\rm Li}_k[x].$$

Immediately from the definitions, we have

\begin{lem}\label{Lem:Stuffle}
\begin{description}
  \item[(1)] For any $w_1,w_2\in\mathfrak{h}^1$, we have
  $${\rm Li}^{\star}[w_1{\bar \ast} w_2]={\rm Li}^{\star}[{w_1}]{\rm Li}^{\star}[w_2].$$
  \item[(2)] Let $n$ be a positive integer and $w_1=\left[\begin{matrix}
k_1\\
x_1
\end{matrix}\right],\ldots,w_n=\left[\begin{matrix}
k_n\\
x_n
\end{matrix}\right],w=\left[\begin{matrix}
k\\
x
\end{matrix}\right]\in\mathcal{M}$. Then we have
$$S\left.\left[\begin{matrix}
k_1, \ldots,k_n\\
x_1,\ldots,x_n
\end{matrix}\right|\begin{matrix}
k\\
x
\end{matrix}\right]={\rm Li}^{\star}[w(w_1{\bar \ast}\cdots{\bar \ast} w_n)].$$
\end{description}
\end{lem}

\pf
One can prove (1) similarly as in \cite{H1997,M2009}, and prove (2) similar as in \cite{Li}.\qed

We prove the corresponding equation of \eqref{Eq:q-Polylog-Product} in the algebra $\mathfrak{h}^1$.

\begin{thm}
Let $n$ be a positive integer and $w_1=\left[\begin{matrix}
k_1\\
x_1
\end{matrix}\right],\ldots,w_n=\left[\begin{matrix}
k_n\\
x_n
\end{matrix}\right]\in\mathcal{M}$. Set $\mathbf{k}=(k_1,\ldots,k_n)$ and $\mathbf{x}=(x_1,\ldots,x_n)$, then we have
\begin{align}
w_1{\bar \ast}\cdots{\bar \ast} w_n=\sum\limits_{\mathbf{l}=(k_{i_1},\ldots,k_{i_m})\in {\rm Sub}(\mathbf{k})\atop \mathbf{l}\neq\mathbf{k}}(-1)^{n-m-1}\left[\begin{matrix}
{\rm wt}(\mathbf{k})-{\rm wt}(\mathbf{l})\\
|\mathbf{x}|/|\mathbf{x}_{\mathbf{l}}|
\end{matrix}\right]\left(w_{i_1}{\bar \ast}\cdots{\bar \ast} w_{i_m}\right).
\label{Eq:Stuffle-Product}
\end{align}
\end{thm}

\pf
We proceed on induction on $n$. The case of $n=1$ is trivial. Now assume that \eqref{Eq:Stuffle-Product} is proved for $\mathbf{k}$ and $\mathbf{x}$. For any $w_{n+1}=\left[\begin{matrix}
k_{n+1}\\
x_{n+1}
\end{matrix}\right]\in\mathcal{M}$, set $\mathbf{k}'=(k_1,\ldots,k_n,k_{n+1})$ and $\mathbf{x}'=(x_1,\ldots,x_n,x_{n+1})$. Using the induction assumption and the definition of the stuffle product, we have
\begin{align*}
&w_1{\bar \ast}\cdots{\bar \ast} w_n{\bar \ast} w_{n+1}\\
=&\sum\limits_{\mathbf{l}=(k_{i_1},\ldots,k_{i_m})\in {\rm Sub}(\mathbf{k})\atop \mathbf{l}\neq\mathbf{k}}(-1)^{n-m-1}w_{n+1}{\bar \ast}\left(\left[\begin{matrix}
{\rm wt}(\mathbf{k})-{\rm wt}(\mathbf{l})\\
|\mathbf{x}|/|\mathbf{x}_{\mathbf{l}}|
\end{matrix}\right]\left(w_{i_1}{\bar \ast}\cdots{\bar \ast} w_{i_m}\right)\right)\\
=&\sum\limits_{\mathbf{l}=(k_{i_1},\ldots,k_{i_m})\in {\rm Sub}(\mathbf{k})\atop \mathbf{l}\neq\mathbf{k}}(-1)^{n-m-1}w_{n+1}\left(\left[\begin{matrix}
{\rm wt}(\mathbf{k})-{\rm wt}(\mathbf{l})\\
|\mathbf{x}|/|\mathbf{x}_{\mathbf{l}}|
\end{matrix}\right]\left(w_{i_1}{\bar \ast}\cdots{\bar \ast} w_{i_m}\right)\right)\\
&+\sum\limits_{\mathbf{l}=(k_{i_1},\ldots,k_{i_m})\in {\rm Sub}(\mathbf{k})\atop \mathbf{l}\neq\mathbf{k}}(-1)^{n-m-1}\left[\begin{matrix}
{\rm wt}(\mathbf{k})-{\rm wt}(\mathbf{l})\\
|\mathbf{x}|/|\mathbf{x}_{\mathbf{l}}|
\end{matrix}\right]\left(w_{n+1}{\bar \ast} w_{i_1}{\bar \ast}\cdots{\bar \ast} w_{i_m}\right)\\
&+\sum\limits_{\mathbf{l}=(k_{i_1},\ldots,k_{i_m})\in {\rm Sub}(\mathbf{k})\atop \mathbf{l}\neq\mathbf{k}}(-1)^{n-m}\left[\begin{matrix}
{\rm wt}(\mathbf{k})-{\rm wt}(\mathbf{l})+k_{n+1}\\
|\mathbf{x}|x/|\mathbf{x}_{\mathbf{l}}|
\end{matrix}\right]\left(w_{i_1}{\bar \ast}\cdots{\bar \ast} w_{i_m}\right)\\
=&w_{n+1}(w_1{\bar \ast}\cdots{\bar \ast} w_n)\\
&+\sum\limits_{\mathbf{l}=(k_{i_1},\ldots,k_{i_m})\in {\rm Sub}(\mathbf{k})\atop \mathbf{l}\neq\mathbf{k}}(-1)^{n-m-1}\left[\begin{matrix}
{\rm wt}(\mathbf{k})-{\rm wt}(\mathbf{l})\\
|\mathbf{x}|/|\mathbf{x}_{\mathbf{l}}|
\end{matrix}\right]\left(w_{n+1}{\bar \ast} w_{i_1}{\bar \ast}\cdots{\bar \ast} w_{i_m}\right)\\
&+\sum\limits_{\mathbf{l}=(k_{i_1},\ldots,k_{i_m})\in {\rm Sub}(\mathbf{k})\atop \mathbf{l}\neq\mathbf{k}}(-1)^{n-m}\left[\begin{matrix}
{\rm wt}(\mathbf{k})-{\rm wt}(\mathbf{l})+k_{n+1}\\
|\mathbf{x}|x/|\mathbf{x}_{\mathbf{l}}|
\end{matrix}\right]\left(w_{i_1}{\bar \ast}\cdots{\bar \ast} w_{i_m}\right).
\end{align*}
Since any $\mathbf{l}\in{\rm Sub}(\mathbf{k}')$ with $\mathbf{l}\neq\mathbf{k}'$ must satisfy and only satisfy one of the following conditions
\begin{description}
  \item[(i)] $\mathbf{l}=(k_1,\ldots,k_n)$;
  \item[(ii)] $\mathbf{l}=(k_{i_1},\ldots,k_{i_m},k_{n+1})$ with $(k_{i_1},\ldots,k_{i_m})\neq\mathbf{k}$;
  \item[(iii)] $\mathbf{l}=(k_{i_1},\ldots,k_{i_m})\neq\mathbf{k}$ and $i_m<n+1$,
\end{description}
we prove \eqref{Eq:Stuffle-Product} for $\mathbf{k}'$ and $\mathbf{x}'$.\qed

Finally, we can prove Theorem \ref{Thm:q-Polylog-Product-New}.

\noindent {\bf Pooof of Theorem \ref{Thm:q-Polylog-Product-New}.}
Applying ${\rm Li}^{\star}$ on both sides of \eqref{Eq:Stuffle-Product}, and with the helps of Lemma \ref{Lem:Stuffle}, we get the result.\qed

\section{Some identities on Euler sums}\label{Sec:Id-EulerSum}

From Theorems \ref{Thm:q-Hurwitz}-\ref{Thm:q-Euler-Second}, taking $x\rightarrow \pm 1,q \rightarrow 1$,  we get the following corollaries.

\begin{cor}(\cite{XC2017})\label{cor4.1}
For positive integers $k>1$ and $l>1$, it holds
\begin{align*}
&(-1)^{k-1}S(1,l;k)-(-1)^{l-1}S(1,k;l)\\
=&\sum\limits_{j=2}^{l-1}(-1)^{j-1}\zeta(l+1-j)S(k;j)-\sum\limits_{j=2}^{k-1}(-1)^{j-1}\zeta(k+1-j)S(l;j) \\
&\quad+\zeta(k)S(1;l)-\zeta(l)S(1;k)+\zeta(l)\zeta(k+1)-\zeta(k)\zeta(l+1).
\end{align*}
\end{cor}

\begin{cor}(\cite{X2017})\label{cor4.2}
For positive integers $k>1$ and $l$, it holds
\begin{align*}
&(-1)^{l-1}S(1,k;l+1)+(-1)^{k-1}S(1,l+1;k)\\
=&\zeta(k+1)\zeta(l+1)+\sum\limits_{j=1}^{k-1}(-1)^{j-1}\zeta(k+1-j)S(j;l+1)+(-1)^{k-1}\zeta(l+1)S(1;k)\\
&\;-\sum\limits_{j=1}^{k-1}(-1)^{j-1}\zeta(k+1-j)\zeta(l+j+1)-\sum\limits_{j=1}^{l-1}(-1)^{j-1}\zeta(l+1-j)S(k;j+1).
\end{align*}
\end{cor}

\begin{cor}\label{cor4.3}
 Let $l \geq 2$ and $k \ge 0$  be integers. Then we have
\begin{align*}
 &(-1)^l\left[S(\bar 1,l+2k+1;l)+S(\bar 1,l;l+2k+1)\right]\\
=&\sum\limits_{j=1}^{l+2k}(-1)^{j-1}\bar\zeta(l+2k+2-j)S(l;\bar j) \\
&\;-\sum\limits_{j=1}^{l-1} (-1)^{j-1}\bar\zeta(l+1-j)S(l+2k+1;\bar j) \\
&\;+(-1)^l\ln 2 \left[S(l+2k+1;l)+S(l;l+2k+1)\right]\\
&\;+(-1)^l\ln 2 \left[S(l+2k+1;\bar l)+S(l;\overline{l+2k+1})\right].
\end{align*}
\end{cor}

\begin{cor}\label{cor4.4}
For integers $l\in \mathbb{N} \setminus \{1\}$ and $k\in \N\cup\{0\}$, we have
\begin{align*}
&(-1)^l\left[S(\bar 1,l+2k;l)-S(\bar 1,l;l+2k)\right] \\
=&\sum\limits_{j=1}^{l+2k-1} (-1)^{j-1}\bar\zeta(l+2k+1-j)S(l;\bar j)  \\
&\;-\sum\limits_{j=1}^{l-1} (-1)^{j-1}\bar\zeta(l+1-j)S(l+2k;\bar j) \\
&\;+(-1)^l\ln 2\left[S(l+2k;l)-S(l;l+2k)\right]\\
&\;+(-1)^l\ln 2\left[S(l+2k;\bar l)-S(l;\overline{l+2k})\right].
\end{align*}
\end{cor}

From Theorem \ref{Thm:q-Polylog-Product}, we find for a positive integer $l>1$, it holds
\begin{align*}
&\zeta^4(l)=4S(\left\{l\right\}_3;l)-6S(\left\{l\right\}_2;2l)+4S(l;3l)-\zeta(4l),\\
&\zeta(2l)\zeta^2(l)=2S(l,2l;l)+S(\left\{l\right\}_2;2l)-S(2l;2l)-2S(l;3l)+\zeta(4l),\\
&\zeta (3l){\zeta ^2}(l) = 2S( {l,3l;l} ) + S( {{{\left\{ l \right\}}_2};3l} ) - S( {3l;2l}) - 2S( {l;4l} ) + \zeta(5l),\\
&\zeta^5(l)=5S(\left\{l\right\}_4;l)-10S(\left\{l\right\}_3;2l)+10S(\left\{l\right\}_2;3l)-5S(l;4l)+\zeta(5l),\\
&\zeta(2l)\zeta^3(l)=S(\left\{l\right\}_3;2l)+3S(\left\{l\right\}_2,2l;l)-3S(\left\{l\right\}_2;3l)-3S(l,2l;2l)\\
&\quad\quad \quad\quad \quad\quad +3S(l;4l)+S(2l;3l)-\zeta(5l).
\end{align*}
Here $\{l\}_d$ denotes the sequence $\underbrace{l,\ldots,l}_{d \text{\;times}}$.

\begin{cor}\label{cor4.5}
For integers $l\in \mathbb{N} \setminus \{1\}$ and $k\in \N\cup \{0\}$, the following identity holds:
\begin{align*}
&S\left( {\bar 1,l + 2k + 1;l} \right) + S\left( {\bar 1,l;l + 2k + 1} \right) + S\left( {l,l + 2k + 1;\bar 1} \right) \\
& =S(l;\overline {l+2k+2})+S(\bar 1;2l+2k+1)+S(l+2k+1;\overline {l+1}) \\
&\quad  +\ln 2\zeta \left( {l + 2k + 1} \right)\zeta \left( l\right) - \bar \zeta \left( {2l + 2k + 2} \right).
\end{align*}
\end{cor}

Hence, from Corollary \ref{cor4.3} and Corollary \ref{cor4.5}, we obtain the following description of quadratic Euler sums.
\begin{cor} \label{cor4.6}
For $l\in \mathbb{N} \setminus \{1\}$ and $k \in \N\cup \{0\}$, the alternating quadratic sums
\[S(l,l + 2k + 1;\bar 1) = \sum\limits_{n = 1}^\infty  {\frac{{H_n^{\left( l \right)}H_n^{\left( {l + 2k + 1} \right)}}}{n}} {\left( { - 1} \right)^{n - 1}}\]
are reducible to linear sums.
\end{cor}
 A simple example is as follows:
\begin{align*}
S(2,3;\bar 1)
 =& - \frac{{161}}{{64}}\zeta( 6 ) + \frac{{31}}{{16}}\zeta ( 5 )\ln 2 + \frac{9}{{32}}{\zeta ^2}( 3) + \frac{3}{8}\zeta( 2 )\zeta ( 3 )\ln 2 + 2\zeta ( 2 ){\rm Li}{_4}\left( {\frac{1}{2}}\right ) \\
&  - \frac{5}{4}\zeta ( 4 ){\ln ^2}2 + \frac{1}{{12}}\zeta ( 2 ){\ln ^4}2 + S(2;\bar 4)  - S(\bar 3;3).
\end{align*}
In fact, proceeding in a similar fashion to evaluation of the Theorem \ref{Thm:q-Euler-First} and Corollary \ref{cor4.6}, it is possible to evaluate other
Euler sums involving harmonic numbers and alternating harmonic numbers. For example, in the same way as in the proof of Corollary \ref{cor4.6}, we also prove that the alternating quadratic sums
\[S(\bar l,\overline {l + 2k + 1};\bar 1) = \sum\limits_{n = 1}^\infty  {\frac{{{\overline H}_n^{\left( l \right)}{\overline H}_n^{\left( {l + 2k + 1} \right)}}}{n}} {\left( { - 1} \right)^{n - 1}}\]
are reducible to linear sums, for $l\in \mathbb{N} \setminus \{1\}$ and $k \in \N\cup \{0\}$. A special case is as follows:
\begin{align*}
S(\bar 2,\bar 3;\bar 1)
& = \frac{{163}}{{128}}\zeta( 6 ) - \frac{{31}}{{16}}\zeta( 5 )\ln 2 + \frac{3}{{16}}{\zeta ^2}( 3 ) - \frac{3}{4}\zeta ( 2 )\zeta ( 3)\ln 2 - \zeta ( 2 ){\rm Li}{_4}\left( {\frac{1}{2}}\right)\\
&\quad + \frac{5}{8}\zeta ( 4 ){\ln ^2}2 - \frac{1}{{24}}\zeta ( 2 ){\ln ^4}2 + S(\bar 2;4) + S(\bar 3;3).
\end{align*}

\noindent {\bf Acknowledgments.} The first author is supported by the National Natural Science Foundation of
China (Grant No. 11471245) and the Natural Science Foundation of Shanghai (grant no. 14ZR1443500). We thank the anonymous referee for suggestions which led to improvements in the exposition.


\end{document}